\documentclass[11pt]{article}
\usepackage[dvips]{epsfig,graphicx}
\usepackage{amsmath,amsfonts,latexsym,enumerate,color}

\font\gotica eufm10  %scaled \magstep1
\font\smallgotica eufm7
\newcommand{\gX}{\mbox{{\gotica X}}}
\newcommand{\sgX}{\mbox{{\smallgotica X}}}
\newcommand{\gU}{\mbox{{\gotica U}}}

\newtheorem{theorem}{Theorem}[section]
\newtheorem{proposition}[theorem]{Proposition}
\newtheorem{lemma}[theorem]{Lemma}
\newtheorem{corollary}[theorem]{Corollary}

\newtheorem{problem}[theorem]{Problem}

\newtheorem{example}[theorem]{Example}

\newcommand{\Qed}{\hfill \rule{2.5mm}{3mm}\medskip}
\newcommand{\proof}{{\sc \noindent Proof.}\ }
\newcommand{\B}{\mathcal{B}}

\newcommand{\Ga}{\Gamma}
\newcommand{\lm}{\lambda}

\newcommand{\inv}{^{-1}}
\renewcommand{\mod}{\hbox{mod}\, }

\newcommand{\ZZ}{\mathbb{Z}}
\newcommand{\QQ}{\mathbb{Q}}

\newcommand{\CC}{\mathbb{C}}
\newcommand{\cay}{{\rm Cay}}

\newcommand{\un}{\underline}

\newcommand{\la}{\langle}
\newcommand{\ra}{\rangle}
\newcommand{\im}{\mathbf{i}}
\newcommand{\diag}{\hbox{{\rm diag}}}
\newcommand{\Det}{\hbox{{\rm det}}}

\newcommand{\GP}{\hbox{{\rm GP}}}

\renewcommand{\dim}{\hbox{{\rm dim}}}
\newcommand{\Ker}{\hbox{{\rm Ker}}}
\newcommand{\Aut}{\hbox{{\rm Aut}}}

\renewcommand{\S}{\mathbf{S}}
\newcommand{\Mat}{\hbox{{\rm Mat}}}
\newcommand{\val}{\hbox{{\rm val}}}

\newcommand{\id}{\hbox{{\rm id}}}

\newcommand{\Sub}{\hbox{{\rm Sub}}}

\newcommand{\comment}[1]{}

\begin{document}
\begin{center}
{\bf\large TRANSITIVE GROUP ACTIONS: (IM)PRIMITIVITY AND SEMIREGULAR SUBGROUPS}
\end{center}
\bigskip

\begin{center}
{\small Istv\'{a}n
Kov\'{a}cs$\,^{a,}$\footnotemark, %\hspace{2pt}\footnotemark,
 ~Aleksander Malni\v c$\,^{a, \, b,}$\addtocounter{footnote}{-1}\footnotemark,
 ~Dragan Maru\v si\v c$\,^{a,}$\hspace{-3pt}\addtocounter{footnote}{-1} \footnotemark,
 \v Stefko Miklavi\v c$\,^{a,}$\hspace{-3pt}\addtocounter{footnote}{-1} \footnotemark $^,$*
 }
 \end{center}

\medskip
\begin{center}
{\em \small
$^a$University of Primorska, Andrej Maru\v{s}i\v{c} Institute, Muzejski trg 2, 6000 Koper, Slovenia and \\
$^b$University of Ljubljana, Pedago\v{s}ka fakulteta, Kardeljeva plo\v{s}\v{c}ad 16, 1000 Ljubljana, Slovenia
}
\end{center}

\addtocounter{footnote}{0}
\footnotetext{Supported in part by ``Javna agencija za raziskovalno dejavnost Republike Slovenije'',
program no. P1-0285.
%\addtocounter{footnote}{1}
%\footnotetext{Supported in part by  the grant OTKA T-043758.

  ~* Corresponding author e-mail: ~stefko.miklavic@upr.si}

\begin{abstract}
The following problem is considered: if $H$ is a semiregular
abelian subgroup of a transitive permutation group $G$ acting on a
finite set $X$, find conditions for (non) existence of
$G$-invariant partitions of $X$.   Conditions  presented in
this paper are derived by studying   spectral properties of
associated $G$-invariant digraphs.  As an essential tool, irreducible
complex characters of $H$  are used.

Questions of this kind  arise naturally when
classifying combinatorial objects which enjoy a certain degree of
symmetry. 
As an illustration, a new
and short proof of an old result of Frucht, Graver and Watkins
({\it Proc. Camb. Phil. Soc.}, {\bf 70} (1971), 211-218)
classifying edge-transitive generalized Petersen graphs, is given.
\end{abstract}

\section{Introductory and historic remarks}
\label{sec:one}
A combinatorial   approach can often provide useful insight in the
study of various algebraic structures, for example,  group actions and
in particular permutation groups. 
On the other hand, when investigating structural  properties of graphs admitting  (transitive)
group actions   one often relies on results
which are purely algebraic. Literature covering both aspects  of  this 
fruitful  interplay is vast, see for example 
\cite{seven, ED06, DMMN, DMMN07, FKW, FGW, GX07, LLZ06,
MNS, DM03, P03a} and the literature therein.

A particular case of such an interplay of  algebraic and combinatorial  concepts
is dealt with in this paper,  and is motivated by the  concept of $B$-groups. 
Recall that a  (finite)group $H$ is called a {\em B-group}
if a primitive permutation group containing
$H$ as a regular subgroup is necessarily doubly transitive.
(A transitive  permutation group $G$ is {\em primitive}  if the two trivial partitions are 
the only block systems of $G$, here called  $G$-invariant partitions; 
$G$  is {\em doubly transitive} provided  it is transitive on 
ordered pairs of distinct elements from the underlying set.)
As a generalization let us call $H$ an \emph{$m$-B-group} if a primitive permutation group
containing $H$ as a semiregular subgroup with $m$ orbits is necessarily doubly transitive.
(A nontrivial permutation, and more generally a
nontrivial permutation group,
is {\em semiregular} if all of its orbits have equal size.)
Of course, a $1$-B-group is just a B-group. By the classical results
of Schur  and Wielandt \cite[Theorem~25.3 and 25.6]{W},
cyclic groups of composite order and dihedral groups are $B$-groups.
Of more recent results we mention the classification of all abelian $B$-groups  by Li \cite{L03},
and  the classification of all metacyclic groups of order $pq$, where $p$ and $q$ are distinct primes, by Dobson \cite{ED06}.

As for $2$-B-groups, cyclic groups of prime order $p$ are of this kind. Namely,  in view of the classification of finite
simple groups,  a primitive group of degree
$2p$, $p>5$ a prime, is necessarily doubly transitive \cite{PC81}. 
We mention that primitive permutations groups containing an element 
with exactly two cycles (not necessarily of the same length) have been classified by M\"uller in \cite{M01}.

As a natural generalization of the original question of Burnside we therefore 
propose the following problem.

\begin{problem}
\label{que:mBgroup}
Given a positive integer $m$, determine the class of (abelian) $m$-B-groups.
\end{problem}

\noindent
We remark  that  asking a transitive permutation group
to contain a semiregular abelian subgroup as in Problem~\ref{que:mBgroup} 
is not that restrictive at all. Indeed, 
although not every transitive permutation group
contains a semiregular subgroup (see for example \cite{seven}), it is believed that the
automorphism group of a  vertex-transitive digraph   does contain such a subgroup
(see \cite[Problem 2.4]{DM81}) and, more generally, that the same holds for 
every $2$-closed transitive permutation group  \cite{bcc15}.
Substantial evidence supporting these claims has been gathered thus far, see \cite{seven, DMMN, DMMN07, MG03, MG07, GX07, KS10, MS98}. 
 
In this paper we develop a general framework within  which the study  
of $m$-B-groups becomes easier to grasp. 
Typically, a natural environment  for building such a framework is graph-theoretic:  the structure of a  transitive
permutation group $G$  containing a semiregular abelian subgroup $H$  is studied via a  collection of  
 $G$-invariant digraphs (orbital graphs)  associated with  $G$.  
 (Recall that a {\em $G$-invariant digraph on $X$}  is
a digraph $\Gamma$ with vertex set $X$ and arc  set a subset of $X \times X$,  
admitting $G$ as subgroup of automorphisms.)
Within  this   setting we 
discuss conditions  which  narrow down 
the possible candidates for nontrivial $G$-invariant partitions, 
thus  paving  the way for  a more in-depth analysis  of $m$-B-groups for small values of 
$m$,  and possible  classification of such groups at least for $m =2$.

From a graph-theoretic point of view  this  line of research is motivated by 
the following situation, often encountered in  problems  regarding 
symmetry of graphs. 
Given a {\em bicirculant}, that is, a graph admitting a cyclic subgroup
with two orbits of equal size, there are two essentially
different possibilities forcing such a graph to be vertex-transitive.
Either there exists
an automorphism swapping setwise the two orbits, called hereafter a {\em swap},
or there exists an automorphism mixing the vertices of the two orbits, called hereafter  a {\em mixer}.
For example, the class of generalized Petersen graphs $GP(n,s)$
\cite{FGW}
contains examples of each of the four possibilities that
may occur with this respect:
first, $GP(7,2)$ is not vertex-transitive and hence without swaps
and mixers, second, the prism $GP(3,1)$ has swaps, but no mixers,
third, the dodecahedron $GP(10,2)$ has mixers, but no swaps,
and finally, the cube $GP(4,1)$ has both swaps and mixers.
It seems therefore natural to seek
for (non)existence conditions for these two kinds of automorphisms
in bicirculants.
Of course, the existence of a mixer in a bicirculant
means that the two orbits of the semiregular cyclic subgroup do not form an invariant partition
of the full automorphism group.
On the other hand, the existence of a swap in a bicirculant is equivalent to there being
a transitive subgroup of automorphisms with the two orbits
of the cyclic group in question forming an invariant partition for this group.
Clearly, knowing which cyclic groups are $2$-B-groups is of vital importance
in the study of arc-transitive bicirculants. On the other hand, 
studying existence of mixers for bicirculants
is essentially the first step in determining which cyclic groups
are $2$-B-groups and which are not.

Coming back to conditions 
for  (non)existence of  certain $G$-invariant partitions,
the  approach we use in this paper is to consider the group $G$ acting on
eigenspaces of $G$-invariant digraphs. 
In order to formulate our main result, some remarks are in order.
First, recall that in a vertex-transitive $G$-invariant digraph  each 
vertex has the same number of in- and out-neighbors, to which
we shall refer as the valency  $\val(\Ga)$  of $\Ga$. 
Further, the partition of $X$ into $H$-orbits  is an
equitable partition for $\Ga$, that is, given two  orbits
$X_1$ and $X_2$, the number of arcs  in $\Ga$ starting from a fixed vertex in
$X_1$ and ending in a vertex in $X_2$ does not depend on the choice of the vertex in $X_1$ (see \cite[p. 76]{G} for details).
Denote this number by $m_{X_1,X_2}$.  The quotient $\Ga/H$ is the  (multi)digraph with vertex
set consisting of all $H$-orbits, and with exactly $m_{X_1,X_2}$ arcs from  $X_1$ to $X_2$. Finally,
the characteristic polynomial of $\Ga/H$ divides the characteristic polynomial of $\Ga$, see \cite[Lemma 2.2(c)]{G}.
A simplified version of our main result can now be stated as follows.
%}

\begin{theorem}
\label{thm:Extreme}
Let $G$ be a transitive permutation group on a finite set $X$
containing an abelian semiregular subgroup $H$. Further, let
$\Gamma$ be a $G$-invariant digraph with
%symbol $\S$ and with
an eigenvalue $\lm$ which is not the valency of $\Ga$.
If the multiplicity of $\lm$ as the eigenvalue of the quotient graph
$\Ga / H$ is the same as the multiplicity of $\lm$ as the eigenvalue of $\Ga$,
then there exists an intransitive normal subgroup of $G$ containing $H$.
\comment{
\textcolor{blue}{
In particular, if $H$ has prime number of orbits, then the orbits of $H$ form
an imprimitivity block system of $G$.}}
\end{theorem}

Observe that, in the particular case when $H$ has  prime number
of orbits, Theorem~\ref{thm:Extreme} gives that the orbits of $H$  coincide with the orbits of 
such an intransitive normal subgroup, and hence 
form a $G$-invariant partition.

\medskip
Theorem~\ref{thm:Extreme} is proved in Section~\ref{sec:five} (as a consequence of more general 
statement, see  Theorem~\ref{thm:moregeneral}), after introducing the
required machinery in Sections~\ref{sec:two} -- \ref{sec:four}.
As an application we give  in Section~\ref{sec:seven}  an alternative proof of an old result
classifying edge-transitive
generalized Petersen graphs \cite{FGW} (see Theorem~\ref{t:gp}).
Theorems of this  nature are  usually proved
either by  an elementary, although technically rather involved,
combinatorial  approach  \cite{FGW,DM98,OW04},
or else by using  a normal subgroup  reduction approach,
often involving the classification of finite simple groups
as an  essential ingredient, see for example
\cite{FP1,FKW,L03,LLZ06,DM03,P99}.
Our proof of Theorem~\ref{t:gp} uses the  second approach. However,  being  based on the
techniques of this paper, it does not rely on the classification of finite simple groups.

%%%%%%%%%%%%%%%%%%%%%%%%%%%%%
%%
%%  SECTION 2
%%
%%%%%%%%%%%%%%%%%%%%%%%%%%%%%

\section{$\boldsymbol{G}$-invariant partitions}
\label{sec:two}

Let $G$ be a permutation group acting trasitively on a finite set $X$.
Suppose that $G$ contains a semiregular subgroup $H$ with orbits $X_i$, $i \in \{1,\dots,m\}$. 
In the next proposition we are going to show that each $G$-invariant partition of $X$ has a particular 
nice structure. To state our result we need the following definition. 
Let $\Delta$ be a partition of $\{1,\dots,m\}$,
and let $\phi \colon \{1,2,\ldots, m\} \to \Sub(H)$, where $\Sub(H)$ denotes the set of subgroups
of $H$. The mapping $\phi$ is {\em $\Delta$-compatible} if $\phi(i) = \phi(j)$ whenever
$i$ and $j$ belong to the same class of $\Delta$.

\begin{proposition}
\label{lem:gtriple}
Let $G$ be a transitive permutation group on a finite set $X$, and let $H$ be a
semiregular abelian subgroup of $G$ with orbits $X_i$, $i \in \{1,\dots,m\}$.
Let $\B$ be a $G$-invariant partition of $X$.
Then there exist some $\un{x}  \in X_1 \times \cdots \times X_m$, a
partition $\Delta$ of $\{1,\ldots,m\}$, and a $\Delta$-compatible mapping
$\phi \colon \{1,\ldots,m\} \to \Sub(H)$ such that
\begin{equation*}
\B = \{\bigcup_{i\in T} x_i^{h\phi(i)}\ |\ T\in  \Delta, \ h \in H\}.
\end{equation*}
\end{proposition}

\proof
Define a relation $\sim$ on $\{1,\dots,m\}$ by letting
$i \sim j$ if and only if there exists some $B \in \B$ such that
$B \cap X_i \ne \emptyset$ and $B \cap X_j \ne \emptyset$.
This is an equivalence relation. It is clearly reflexive and symmetric.
To see that it is also transitive,  let $i \sim j$ and $j \sim k$. 
Then there are  $B_1,B_2 \in \B$  such that $B_1 \cap X_i \ne \emptyset$,
$B_1 \cap X_j \ne \emptyset$ and 
$B_2 \cap X_j \ne \emptyset$, $B_2 \cap X_k \ne \emptyset$. 
Let $x \in B_1 \cap X_j$ and $x'  \in B_2 \cap X_j$.
Since $X_j$ is an orbit of $H$, there is
$h\in H$ such that $x'=x^h$.
As $\B$ is a $G$-invariant partition of $X$ we have  $B_1^h=B_2$.
But $X_i$ is also an orbit for $H$; hence
 $B_1 \cap X_i \ne \emptyset $ implies
 $B_1^h \cap X_i^h \ne \emptyset $,
that is, $B_2 \cap X_i \ne \emptyset $.
Since $B_2 \cap X_k \ne \emptyset$, we have $i \sim k$. 

Let $\Delta$ denote the partition of $\{1,\dots,m\}$
induced by the equivalence classes of $\sim$. 
Note that in the proof of transitivity of $\sim$ we actually prove a
stronger property: if for some $B \in \B$ and some $T \in \Delta$
we have $B \cap X_i \ne \emptyset$ for some $i \in T$, then
$B \cap X_i \ne \emptyset$ for every $i \in T$.  
This allows us to make the following definition.
For each  $T\in \Delta$  set
$$
  \B_T=\{\,  B\in \B \mid B \cap X_i \neq \emptyset \mbox{ for all }i\in T\, \},
$$
fix  a block  $B_T \in \B_T$, and for each $i \in T$ choose an element $x_i \in B_T \cap X_i$.
Denoting by  $K_T$ the setwise stabilizer of $B_T$ in $H$, let
$\phi \colon \{1, \ldots, m\} \to \Sub(H)$ be the mapping defined by
$$%
\phi(i) =K_T,  \quad  T \in \Delta, \ i \in T.
$$%
We now show that 
$$
 \B=\{\bigcup_{i\in T} x_i^{h\phi(i)}\ |\ T\in  \Delta, \ h \in H\}.
$$
First, observe that for each $i \in T$ we have
$$
K_T=\{ h \in H \mid x_i^h \in B_T \}.
$$
Next, let  $B'$ be an arbitrary block in $\B_T$ and let $i \in T$.
Since $X_i$ is an $H$-orbit  there exists $h\in H$ such that $B'=B_T^h$. Therefore,
as  $H$ is abelian and $B_T = \cup_{i \in T} x_i^{K_T}$ we have
$B' = \bigcup_{i\in T} x_i^{hK_T}$.
This completes the proof.
\Qed

For later reference we shall denote the description of $G$-invariant partition $\B$ 
in Proposition~\ref{lem:gtriple} by $\Pi(\un{x}, \Delta, \phi)$.

\begin{corollary}
\label{cor:gtriple}
With notation and assumptions of Proposition \ref{lem:gtriple},  let the partition 
$\B = \Pi(\underline{x}, \Delta, \phi)$
be $G$-invariant. Then the following hold.
\begin{itemize}
\item[(i)] 
If  $N$ is the kernel of the action of $G$ on $\B$, then  for each  $j \in \{1, \ldots, m\}$
we have $H \cap N \le \phi(j)$.
\item[(ii)] 
For $T, T' \in \Delta$ and $j \in T$, $j' \in T'$ we have that if $|\phi(j)|=|\phi(j')|$  then  $|T|=|T'|$ .
\end{itemize}
\end{corollary}

\proof 
To prove (i), pick $j \in \{1,2, \ldots, m\}$ and $T \in \Delta$ such that $j \in T$. 
Observe that if $h \in H \cap N$, then $x_j^h \in B_T$. 
This implies $h \in \phi(j)$, that is, $H \cap N \le \phi(j)$.

To prove (ii),  note that since all blocks of $\B$ are of the same size, we have
$|\phi(j)| |T|=|\phi(j')| |T'|$ for all $T, T' \in \Delta$ and $j \in T$, $j' \in T'$.
The result now follows. \Qed

%%%%%%%%%%%%%%%%%%%%%
%
% SECTION 3
%
%%%%%%%%%%%%%%%%%%%%

\section{ Eigenspaces of $\boldsymbol{G}$-invariant digraphs }
\label{sec:three}

Let $G$ be a transitive permutation group on a finite set $X$.
An orbit of $G$ acting on $X \times X$ is called an {\em orbital} of $G$.
A {\em $G$-invariant digraph} over $X$ is a digraph with $X$ as the set of vertices 
such that $G$ acts as a group of automorphisms of $X$. The set of arcs of a  
$G$-invariant digraph is necessarily a union of orbitals.
In this section we describe the eigenspaces of $G$-invariant digraphs in the case when $G$ has
an abelian semiregular subgroup. The results of this section will be used in Section~\ref{sec:four} 
where we analyze $G$-invariant partitions of $X$.

\medskip
We start by  introducing  some notation. Let $Z$ be  a
nonempty finite set  with a fixed linear ordering of
its elements. By  $V_Z$ we denote the $|Z|$-dimensional $\CC$-vector space of row
vectors over $\CC$, with components of vectors  indexed  by the elements of $Z$ 
(relative to the chosen ordering of $Z$).
In particular, we  occasionally consider  a complex valued function
$f \colon Z \to \CC$ as a vector  $\mathbf{v}_f \in V_Z$ whose
$z$-th component,  $z \in Z$, is equal to $f(z)$.
Similarly, by  $\Mat_Z$ we  denote the algebra of $|Z| \times |Z|$ matrices, with
columns and rows indexed by the elements of $Z$ (relative to the chosen ordering of $Z$).

Recall that for a subset $S$ of a group $H$, the {\em Cayley digraph} $\cay(H,S)$
is the digraph with vertex set $H$, and with  an arc
pointing from $x\in H$ to $y\in H$ whenever $yx\inv \in S$.
Fixing   a linear ordering  $h_1,h_2, \dots ,h_n$ of  elements of $H$ we let
$A(\cay(H,S)) \in \Mat_H$  denote
the {\em adjacency matrix}  of $\cay(H,S)$  whose  $(h_i, h_j)$-th entry
is $1$ if $h_jh_i\inv \in S$, and is $0$ otherwise. If $H$ is abelian, let $H^*$
denote the {\em dual group} of complex irreducible characters  $\chi\colon H \to \CC$. 
By \cite[Lemma~9.2]{G} and Remark below, $\mathbf{v}_\chi$ is an eigenvector of $A(\cay(H,S))$; 
its corresponding eigenvalue is 
$$%
\chi(S)=  \sum_{g\in S}\chi(g^{-1}) = \sum_{g\in S} \overline{\chi(g)},
$$%
with the standard convention that  $\chi(S) = 0$ whenever $S = \emptyset$.
Moreover, the spectrum of  $A(\cay(H,S))$  is the multiset
$\{\chi(S)\ |\ \chi \in H^*\}$,  and $\{\mathbf{v}_\chi$,  $\chi \in H^*\}$
forms  an orthogonal  basis  of $V_H$ relative to the dot product.
Note that this basis does not depend on $S$.

\medskip
\noindent
{\bf Remark}. Since all  group actions in this paper are on the right we need to consider matrices acting 
on row vectors. It is for this reason that the eigenvalues of $\cay(H,S)$ are as  defined above. 
However, in the  context of left  actions on column vectors, see for example \cite[Lemma~9.2]{G},
the eigenvalues  are given by $\sum_{g\in S}\chi(g)$.

\medskip
The above results can be stated in a more general setting.
Let $G$ be a transitive permutation group on a finite set $X$ and let $H \leq G$
be a semiregular subgroup. Choose a  fixed linear ordering $h_1,h_2, \dots ,h_n$ of elements
of $H$ and a linear ordering  $X_1,\dots, X_m$ of the set $\gX$ of  $H$-orbits on $X$, along with
arbitrarily chosen {\em points of reference}  $b_i \in X_i$ for each $i$.
This naturally induces a linear ordering of $X$, namely
$$%
b_1^{h_1},\dots,b_1^{h_n}, b_2^{h_1},\dots,b_2^{h_n},\dots, b_m^{h_1},\dots,b_m^{h_n}.
$$%

Now let $\Ga$ be a $G$-invariant digraph on $X$.
For $x, y \in X$ we shall write $x \to_\Ga y$ if there is an arc pointing from $x$ to $y$.
For $1 \le i,j \le m$ define
$$
 S_{ij} =  \{h \in H \mid b_i \to_\Ga b_j^h \}.
$$
The array $\S= [S_{ij} \mid 1 \le i,j \le m]$ is called the
{\em symbol of} $\Ga$ relative  to  $H$, the ordering of $\gX$, and the chosen points of
reference  $(b_1, b_2, \ldots, b_m)$.
The reader may check that the  adjacency matrix
of $\Ga$ (relative to the chosen orderings)  is given by the  block matrix

$$%
  A_\S=\begin{pmatrix} A(\cay(H,S_{11})) & \cdots & A(\cay(H,S_{1m})) \\
                       \vdots            & \ddots & \vdots \\
                       A(\cay(H,S_{m1})) & \cdots & A(\cay(H,S_{mm}))
       \end{pmatrix}.
$$%
In order to display $A_\S$ in a more compact form,  let $E_{ij}$ be the elementary 
$m \times m$ matrix 
having the $(i,j)$-th entry equal to $1$ and all other entries equal to $0$, and let 
$A_{ij} = A(\cay(H,S_{ij})$. Then 
\begin{equation*}
A_\S = \sum_{i,j = 1}^m E_{ij} \otimes A_{ij},
\end{equation*}
where $\otimes$ denotes the  tensor (Kronecker)  product of matrices. By \cite[p.~107]{Eves},
$(A \otimes B) (C \otimes D) = AC \otimes BD$;  in particular,   
$A^{-1} \otimes B^{-1} =  (A \otimes B)^{-1}$. 
For later use,  observe that for   any  $x=b_i^{h_j} \in X$, the  $x$-th component of the tensor 
product  $\mathbf{u} \otimes \mathbf{v}_f \in V_X$,  where $\mathbf{u} = (u_1,\dots,u_m) \in \CC^m$ and 
$f\colon H \to  \CC$,    is  given  by
$$
(\mathbf{u} \otimes \mathbf{v}_f)_x=u_if(h_j).  
$$
Further, if $U$ and $V$ are subspaces of $\CC^m$ and $V_H$, respectively, then
$U \otimes V$ is  the subspace of $V_X$ spanned by
all $\mathbf{u} \otimes \mathbf{v}$, $\mathbf{u} \in U$ and $\mathbf{v} \in V$.
Finally, for a   character $\chi \in H^*$ we define the
$m \times m$ matrix $\chi(\S) \in \Mat_{\sgX}$  by
$$%
  \chi(\S)=\begin{pmatrix} \chi(S_{11}) & \cdots & \chi(S_{1m}) \\
                           \vdots       & \ddots & \vdots       \\
                           \chi(S_{m1}) & \cdots & \chi(S_{mm}) \end{pmatrix}.
$$%
Note that  $\chi_0(\S)$, where $\chi_0$ is the principal character of
$H$,  gives us   the  adjacency matrix of the quotient (multi)digraph $\Ga/H$ 
relative to the chosen ordering of  $\gX$.

\begin{proposition}
\label{cor:spec}
With notation and assumptions above, 
let $G$ be a transitive permutation group on a finite set $X$ and
let $H \leq G$ be a semiregular abelian subgroup of  order $n$  with $m$ orbits on $X$. Then the spectrum of a
$G$-invariant digraph $\Ga$ on $X$ is equal to the union (counting multiplicities) of spectra of all
$\chi(\S)$, $\chi \in H^*$,  where  $\S$  is  the  respective
symbol of $\Ga$.
\end{proposition}

\proof
The spectrum of $\Ga$ is equal to the spectrum of
$A_\S$.  Recall that   all  blocks $A_{ij}=A(\cay(H,S_{ij}))$
have  a common orthogonal eigenvector basis  $\{\mathbf{v}_\chi\  | \ \chi \in H^*\}$,
and that the  eigenvalue of $A_{ij}$ corresponding to $\mathbf{v}_{\chi}$ is equal to $\chi(S_{ij})$.

Let us fix an ordering $\chi_1,\ldots,\chi_n$ of $H^*$.
With respect to  the ordered  eigenbasis $\mathbf{v}_{\chi_1}, \ldots, \mathbf{v}_{\chi_n}$, the matrix $A_{ij}$
takes the  diagonal form $D_{ij} = \diag(\chi_1(S_{ij}), \dots, \chi_n(S_{ij}))$, that is,
$A_{ij} = P^{-1} D_{ij} P$ where  $P$ is the $n \times n$ matrix with
$\mathbf{v}_{\chi_1}, \ldots, \mathbf{v}_{\chi_n}$ as rows.  Consequently,
$$%
A_\S = (I_m \otimes P)^{-1} (\sum_{i,j=1}^{m} E_{ij} \otimes D_{ij}) (I_m \otimes P),
$$%
where  the matrix  $I_m$ is the identity matrix of size $m$.
 The reader may check that there is a permutation matrix which establishes a similarity relation
between
$$%
\sum_{i,j=1}^{m} E_{ij} \otimes D_{ij} \quad \hbox{{\rm  and}} \quad \sum_{k = 1}^n E_{kk} \otimes \chi_k(\S).
$$%
As the latter of the two  matrices  is block-diagonal, it is now obvious that  the spectrum of $A_\S$ is equal 
to the union of spectra of all matrices $\chi_k(\S)$.
\Qed

\bigskip
To calculate  the eigenvectors of $A_\S$,  let  $\lm$ be
an eigenvalue and let  $W_{\lm} \leq V_X$ denote
the respective eigenspace.
First, define  the set
\begin{equation}
\label{eq:defK}
K_{\S,\lm} = \big\{\, \chi \in H^*   \mid    \Det(\, \chi(\S)-\lm I_m\, )=0\, \big\} \nonumber
\end{equation}

\noindent
of all  characters $ \chi \in H^*$ for which $\lm$ is an eigenvalue of $\chi(\S)$.
For $\chi \in K_{\S,\lm}$,  let $V_{\lm,\chi}$ be  the $\lm$-eigenspace of the
matrix $\chi(\S)$. Now the eigenspace  $W_{\lm}$  can be written as  
$W_\lm=\bigoplus_{\chi\in K_{\S,\lm}}\, V_{\lm,\chi} \otimes \la \mathbf{v}_\chi \ra$.
Indeed, choosing an orthogonal basis  
$\gU_{\lm,\chi}$ of $V_{\lm,\chi}$ we have the following proposition.

\begin{proposition}\label{cor:espace}
With notation and assumptions above,
if $\lm$ is an eigenvalue of $\Ga$, 
a basis of $W_{\lm}$  is given by 
$$
\{\mathbf{u} \otimes \mathbf{v}_{\chi}\ |\  \chi \in K_{\S,\lm}, \ \mathbf{u} \in \gU_{\lm,\chi} \}.
$$
\end{proposition}

\proof
 Fix an ordering $\chi_1,\ldots,\chi_n$ of $H^*$, and recall
(from the proof of Proposition~\ref{cor:spec}) that $A_\S$ and
$\sum_{k = 1}^n E_{kk} \otimes \chi_k(\S)$  are similar. Thus, if
 $\lm$ is  an eigenvalue of $A_\S$, the dimension of $W_{\lm}$ is equal to
 the dimension of the corresponding
$\lm$-eigenspace of
$\sum_{k = 1}^n E_{kk} \otimes \chi_k(\S)$. Since this $\lm$-subspace  is equal to
$\oplus_{\chi_k \in K_{\S,\lm}} \la \mathbf{e}_k\ra \otimes V_{\lm,\chi_k}$, where
$\{ \mathbf{e}_1, \ldots, \mathbf{e}_n\}$ is the standard basis of $\CC^n$,
we have
$$%
\dim\,W_{\lm} = \sum_{\chi_k \in K_{\S,\lm}} \dim\,V_{\lm,\chi_k}.
$$
To complete the proof we need to find an appropriate basis of $W_{\lm}$. First note that
for $\chi \in K_{\S,\lm}$ and   $\mathbf{u} \in V_{\lm,\chi}$ the vector
$\mathbf{u} \otimes \mathbf{v}_\chi$ is a $\lm$-eigenvector of $A_\S$. Indeed,
\begin{eqnarray*}
(\mathbf{u} \otimes \mathbf{v}_\chi) A_\S  & = &
(\mathbf{u} \otimes \mathbf{v}_\chi) \bigg(\sum_{i,j=1}^{m} (E_{ij} \otimes A_{ij})\bigg)  \\
& = &  \sum_{i,j=1}^{m} (\mathbf{u}\, E_{ij} \otimes \mathbf{v}_\chi A_{ij} )
 =  \sum_{i,j=1}^{m} \big(\mathbf{u}\, E_{ij}  \otimes \chi(S_{ij}) \mathbf{v}_{\chi}\big) \\
 & = &  \bigg(\mathbf{u} \sum_{i,j=1}^{m}  \chi(S_{ij}) \, E_{ij}\bigg) \otimes \mathbf{v}_\chi
  =   \mathbf{u}\, \chi(\S)  \otimes \mathbf{v}_\chi \\ \\
 & = & \lm (\mathbf{u} \otimes \mathbf{v}_\chi).
\end{eqnarray*}
%For each $\chi \in K_{\S,\lm}$ choose an orthogonal basis
%$\gU_{\lm,\chi}$ of  $V_{\lm,\chi}$.
%= \{\mathbf{u}_{\lm,1}^{(k)}, \ldots, \mathbf{u}_{\lm,r_k}^{(k)}\}$.
As the set
$\{\mathbf{u} \otimes \mathbf{v}_{\chi}\ |\  \chi \in K_{\S,\lm}, \ \mathbf{u} \in \gU_{\lm,\chi} \}$
contains  exactly $\sum_{\chi \in K_{\S,\lm}} \dim\,V_{\lm,\chi}$
elements, it remains to show that  it is linearly independent.
To this end we show that the vectors of this set are pairwise orthogonal.
Consider two distinct eigenvectors $\mathbf{u} \otimes \mathbf{v}_{\chi}$ and
$\mathbf{u'} \otimes \mathbf{v}_{\chi'}$ where $\mathbf{u} \in \gU_{\lm,\chi}$ and
$\mathbf{u'} \in \gU_{\lm,\chi'}$.
Their dot  product is
$$%
\la \mathbf{u} \otimes \mathbf{v}_{\chi}, \mathbf{u'} \otimes \mathbf{v}_{\chi'} \ra =
\la \mathbf{u}, \mathbf{u'} \ra \la  \mathbf{v}_{\chi}, \mathbf{v}_{\chi'} \ra.
$$%
If $\chi \neq \chi'$, then  $\la  \mathbf{v}_{\chi}, \mathbf{v}_{\chi'} \ra = 0$.
If $\chi = \chi'$ and $\mathbf{u} \neq \mathbf{u'}$, then $\la \mathbf{u}, \mathbf{u'} \ra = 0$ since
$\gU_{\lm,\chi}$ is an orthogonal basis for $V_{\lambda, \chi}$.
Therefore $\la \mathbf{u} \otimes \mathbf{v}_{\chi}, \mathbf{u'} \otimes \mathbf{v}_{\chi'} \ra = 0$.
The proof is complete.
\Qed

\begin{example}
\label{ex:cube0}
{\rm
Let $\Gamma = Q_3$ be  the $3$-cube,   and let  $c \in \Aut(Q_3)$ be the central reflection
switching  antipodal vertices. Let 
 $G \le \Aut(Q_3)$  be a  vertex-transitive
subgroup containing $c$. 
For the semiregular subgroup $H \leq G$ we take $H = \{\id, c\}$. 
Let $X_1$, $X_2$, $X_3$ and $X_4$ be the vertex orbits of $H$.
Identifying $H  = \ZZ_2$   and choosing  the points of reference 
$b_j \in X_j$  $(j = 1,2,3,4)$ in such a way that no two of them are  adjacent,  the  respective symbol is  given by 
$S_{jj} = \emptyset$ and $S_{ij} = \{1\}$,  for $i \neq j$. 
The elements of the dual group $H^*$ are  $\chi_0(j)= 1$ and $\chi_1(j) = (-1)^j$, for $j \in \ZZ_2$.

By direct computation we get that the eigenvalues of  
$\chi_0(\S)$ are $3, -1,-1,-1$, and  the eigenvalues of $\chi_1(\S)$ are $-3,1,1,1$.
Therefore, 
$$%
K_{\S,3} =  K_{\S,-1} = \{\chi_0\}  \quad \text{and} \quad  \ K_{\S,1}= K_{\S,-3} = \{\chi_1\}.   
$$%
The respective subspaces
$V_{\lm,\chi}$  are
\begin{eqnarray*}
               V_{3, \chi_0} = V_{-3,\chi_1} & = & \la\, (1,1,1,1) \,\ra  \\
               V_{1, \chi_1} = V_{-1,\chi_0} & = & \la\, (1,0,0,-1), (0,1, 0,-1),  (0,0,1,-1) \,\ra.
 \end{eqnarray*}              
               
\noindent            
By Propositions~\ref{cor:spec} and \ref{cor:espace},  the spectrum of $Q_3$ is
$\{3, 1,1,1, -1,-1,-1, -3\}$, and the eigenspaces $W_\lm$ are:

\begin{eqnarray*}
W_3 & = & \la\;(1,1,1,1,1,1,1,1)\,\ra \\
W_1 & = & \la\;(1,-1,0,0,0,0,-1,1),  (0,0,1,-1,0,0,-1,1),  (0,0,0,0,1,-1,-1,1)\,\ra  \\
W_{-1} & = & \la\;(1,1,0,0,0,0,-1,-1),  (0,0,1,1,0,0,-1,-1),  (0,0,0,0,1,1,-1,-1)\,\ra  \\
W_{-3} & = & \la\;(1,-1,1,-1, 1,-1,1,-1)\,\ra.  
\end{eqnarray*}
}
\end{example}

%%%%%%%%%%%%%%%%%%
%
% SECTION 4
%
%%%%%%%%%%%%%%%%%%

\section{$\boldsymbol{G}$-invariant partitions from linear actions}
\label{sec:four}

With  the notation of Section~\ref{sec:three} we  extend the action of $G$ on $X$ 
to an action of $G$ on the vector space of row vectors $V_X$ in a natural way by setting
$$
(\mathbf{w}^g)_x=\mathbf{w}_{x^{g^{-1}}}, \; x \in X, \; g \in G.
$$
Note  that $\mathbf{w}^g=\mathbf{w} P_g$, where  $P_g \in \Mat_X$ is the permutation
matrix defined by 
$$
(P_g)_{x,y}=\left\{ \begin{array}{rl} 1, & x^g=y \\
0,&\textrm{otherwise} \end{array}\right.  \qquad x,y \in X.
$$
Let $\lm$ be an eigenvalue of a $G$-invariant digraph $\Gamma$.  Since $P_g$ commutes with the 
adjacency matrix $A_\S$, the eigenspace $W_\lm$ is clearly invariant for the action of $G$ on $V_X$.
Although $G$  acts on $V_X$ faithfully, this might not be the case with $G$ acting on $W_\lm$.
Denote by $N_\lm$ the kernel of the action of $G$ on $W_\lm$.
The orbits of $N_\lm$ clearly form a $G$-invariant partition $\B_{\lm}$ of $X$.
By Proposition~\ref{lem:gtriple},
$$%
\B_{\lm}=\Pi(\un{x},\Delta_\lm,\phi_\lm)
$$%\
for an appropriate choice of $\un{x} = (x_1, x_2, \ldots, x_m)$, $\Delta_\lm$, and $\phi_\lm$.
Our goal in this section is to describe $\B_{\lm}$ in terms of $\Delta_\lm$ and $\phi_\lm$.

\medskip
We first describe the function $\phi_\lm$. To this end we introduce the following notation.  
 For an ordered  basis $\gU_{\lm,\chi}=(\mathbf{u}_1,...,\mathbf{u}_r)$ of $V_{\lm,\chi},$  we let
$\gU_{\lm,\chi}^{(k)}$ denote the vector of length $r$ 
whose $i$-th component is equal to the $k$-th component of $\mathbf{u}_i \in \gU_{\lm,\chi}$.

\begin{lemma}
\label{lem:K}
With notation and assumptions above the following hold.
\begin{enumerate}[(i)]
\item
$N_\lm \cap H = \cap_{\chi \in K_{\S,\lm}} \Ker\chi$.
\item
The kernel of the action of $G$ on $\B_{\lm}$ is equal to $N_\lm$.
\item
Let  $i \in \{1, \ldots, m\}$ and $\chi \in K_{\S,\lm}$. If    $\gU_{\lm,\chi}^{(i)} \neq 0$, then 
$\phi_\lm(i) \leq \Ker\chi$. In particular, if  $\gU_{\lm,\chi}^{(i)} \neq 0$  for all 
$\chi \in K_{\S,\lm}$, then  $\phi_\lm(i)=N_\lm \cap H$.
\end{enumerate}
\end{lemma}

\proof
First  observe that
 $$
(\mathbf{u} \otimes \mathbf{v}_\chi)^h=
\chi(h^{-1})(\mathbf{u} \otimes \mathbf{v}_\chi)
$$
for all  $\chi \in K_{\S,\lm}$, $\mathbf{u} \in V_{\lm,\chi}$, and $h \in H$.
Indeed,  take an arbitrary   $x=b_i^{h'}$, where  $h' \in H$. Then
\begin{eqnarray*}
((\mathbf{u} \otimes \mathbf{v}_\chi)^h)_x & = &  (\mathbf{u} \otimes \mathbf{v}_\chi)_{x^{h^{-1}}}  \ = \ 
(\mathbf{u} \otimes \mathbf{v}_\chi)_{b_i^{h'h^{-1}}} \\
& = & u_i\chi(h'h^{-1}) \ = \  \chi(h^{-1}) \big(u_i \chi(h') \big)  \\
& = &  \chi(h^{-1})(\mathbf{u} \otimes \mathbf{v}_\chi)_x. 
\end{eqnarray*}
By Proposition~\ref{cor:espace} we have  that $h \in H$  belongs to $N_\lm$
if and only if
$(\mathbf{u} \otimes \mathbf{v}_\chi)^h=\mathbf{u} \otimes \mathbf{v}_\chi$
for each $\chi \in K_{\S,\lm}$ and $\mathbf{u} \in V_{\lm,\chi}$.
Hence  $h \in N_\lm$ if and only if
$\overline{\chi(h)} = \chi(h^{-1}) = 1$ for all $\chi \in K_{\S,\lm}$, or equivalently,
$h \in \cap_{\chi \in K_{\S,\lm}} \Ker \chi$. Part (i) follows.

To show (ii),  let $g \in G$ acting on $X$  preserve  the orbits in $\B_\lm$. 
Recall that by the definition of $N_\lm$ the coordinates of $\mathbf{w} \in W_{\lm}$,
indexed by the elements from the same orbit of $N_\lm$, are the same. Thus, 
$\mathbf{w}_{x^{g^{-1}}} = \mathbf{w}_x$ for all $\mathbf{w} \in W_{\lm}$ and $x \in X$,
and so  $(\mathbf{w}^g)_x=\mathbf{w}_{x^{g^{-1}}} = \mathbf{w}_x$.  Consequently, 
$\mathbf{w}^g = \mathbf{w}$ and hence  $g \in N_\lm$.
As  the kernel of the action of $G$ on $\B_{\lm}$ contains $N_\lm$, the result follows.

In order to show (iii), note that since $\gU_{\lm,\chi}^{(i)} \neq 0$, there exists 
$\mathbf{u} \in V_{\lambda,\chi}$ with $u_i \neq 0$. 
Let $B$ denote the block of $\B_{\lm}$ containing $x_i$.
Recall from the proof of Proposition~\ref{lem:gtriple}  that $\phi_{\lm}(i)=\{h \in H \mid x_i^h \in B\}$.
Fix $h \in \phi_{\lm}(i)$.
Since $x_i$ and $x_i^h$ are in the same $N_\lm$-orbit, the
$x_i$-th and $x_i^h$-th coordinate of $\mathbf{u} \otimes \mathbf{v}_{\chi}$ are the
same. As $x_i=b_i^{h'}$ for some $h' \in H$ we find that
$u_i \chi(h') = u_i \chi(h'h) $. Now $u_i \ne 0$ implies $\chi(h)=1$, 
and consequently $\phi_{\lm}(i) \le \Ker \chi$.
Now if $\gU_{\lm,\chi}^{(i)} \neq 0$ for every
$\chi \in K_{\S,\lm}$, we get
$\phi_{\lm}(i) \le \cap_{\chi \in K_{\S,\lm}} \Ker\chi=N_\lm \cap H$, see (i).
But since also $N_\lm \cap H \le \phi_{\lm}(i)$, see  Corollary~\ref{cor:gtriple}(i),
the claim is proved.
\Qed

\noindent
Following \cite{BS}, for 
subsets $L \subseteq H$ and $L^* \subseteq H^*$ we define 
$$
L^\bot = \{\, \chi\in H^* \mid  L \subseteq \Ker\chi\,\}  \quad \hbox{{\rm and}} \quad
L^{*\bot} = \cap_{\chi\in L^*}\Ker \chi.
$$
Observe that $L^\bot = \la L \ra^\bot$ and $L^{*\bot} = \la L^* \ra^\bot$. Moreover,  
$\bot$ is an antiisomorphism between the subgroup lattices of $H$ and $H^*$. 
In this notation, Lemma~\ref{lem:K}(i) can be restated as
$$ N_\lm \cap H=\la K_{\S,\lm} \ra^\bot.$$

\medskip
We now examine the partition $\Delta_\lm$.  Since  $\Delta_\lm$ is not
readily at hand, we introduce  certain auxiliary partitions of  $\{1,\dots,m\}$
which  are  easier to compute, and, as it will
become clear shortly,   they all have  $\Delta_{\lm}$ as a refinement.
Hence their intersection  gives  at least some information
about $\Delta_{\lm}$, the more so when the former  is `close' to the trivial partition 
$\Delta_I$ consisting of singletons, and giving no information if this intersection
is the partition $\Delta_U$ consisting of the whole set.

For $\chi \in K_{\S,\lm}$ choose an ordered basis $\gU = \gU_{\lm,\chi}$ of $V_{\lm,\chi}$. 
Let $\sim_{\lm,\chi}$ be the equivalence relation on $\{1,\dots,m\}$ defined by
$$
i \ \sim_{\lm,\chi} \  j \iff  \gU^{(i)} =  \chi(h)\gU^{(j)} \
\ \hbox{{\rm for some}} \ \ h \in H.
$$
Note that this relation does not depend on the chosen basis $\gU$. 
We define the auxiliary partition $\Delta_{\lm,\chi}$ as the equivalence classes of
the relation $\sim_{\lm,\chi}$.

\begin{lemma}\label{lem:Delta}
With notation and assumptions above, for $\chi \in K_{\S,\lm}$
we have that $\Delta_\lm$ is a refinement of $\Delta_{\lm,\chi}$.
\end{lemma}

\proof
Fix  a  character $\chi \in K_{\S,\lm}$ and an ordered basis $\gU$  of $V_{\lm,\chi}$.
Let $i$ and $j$ be in the same partition class of $\Delta_\lm$.
Then there exists $h \in H$ such that $b_i$ and $b_j^h$ are in the same $N_\lm$-orbit.
Let $\mathbf{u} = (u_1,u_2,\dots, u_m) \in \gU$.
By Proposition~\ref{cor:espace}, $\mathbf{u} \otimes \mathbf{v}_{\chi} \in W_{\lm}$.
Therefore its $b_i$-th and $b_j^h$-th coordinates are the
same, implying $u_i = \chi(h)u_j$. Since  $\mathbf{u} \in \gU$ was arbitrary, it follows that
$\gU^{(i)} =  \chi(h)\gU^{(j)}$. Hence $i$ and $j$ are in the same partition class of 
$\Delta_{\lm,\chi}$, as required.
 \Qed

\bigskip
Let $\B_\lm =\Pi(\un{x},\Delta_\lm,\phi_\lm)$ be the $G$-invariant partition arising as orbits of 
the normal subgroup $N_\lm \leq G$. Recall from the proof of Proposition~\ref{lem:gtriple}
that one cannot determine $\un{x}$ without knowing $\B_\lm$. In contrast to this, 
$\Delta_\lm$ and  $\phi_\lm$  can be computed -- at least in certain cases  --  and  this is sometimes even sufficient in 
order to completely reconstruct $\B_\lm$.  These remarks are illustrated in the next example.

\begin{example}
\label{ex:cube}
{\rm
We  continue with  the analysis of $\Gamma = Q_3$ from Example~\ref{ex:cube0}.
For each eigenvalue $\lm$ of $Q_3$ 
we are going to find the invariant partition $\B_\lm =\Pi(\un{x},\Delta_\lm,\phi_\lm)$ of the vertex set of $Q_3$.   

By Lemma~\ref{lem:K}(iii) we have that $\phi_\lm(i) = N_\lm \cap H$ for each eigenvalue $\lambda$ and 
$i \in \{1,2,3,4\}$. 
Finding $N_\lm \cap H$ is not a problem. 
Using Lemma \ref{lem:K}(i)  we obtain
 $$
 N_3 \cap H= N_{-1} \cap H = H, \ \  N_{-3} \cap H= N_1 \cap H = \{0\}.
$$

In order to find $\Delta_\lm$ 
we first compute the auxiliary partitions $\Delta_{\lm,\chi}$. By definition we easily get that 
$$%
\Delta_{1,\chi_1} = \Delta_{-1,\chi_0} = \Delta_I, \quad 
\Delta_{3,\chi_0} = \Delta_{-3,\chi_1} = \Delta_U.
$$%
Therefore, by Lemma~\ref{lem:Delta}  we obtain 
that $\Delta_1 = \Delta_{-1} = \Delta_I$.  Consequently, $N_1 \cap H = \{0\}$ implies that 
$\B_1$ is the trivial block system consisting of singletons, and  $N_{-1} \cap H = H$ implies that
the blocks of $\B_{-1}$ are the orbits of $H$.  

In the remaining two cases,  direct analysis of the action of the 
group $G$ on $W_{\lm}$ is needed in order to compute  $N_{\lm}$, and hence $\B_\lm$.
Consider $W_3=\langle (1,1,1,1,1,1,1,1)\rangle$. Then, obviously,  
$N_3=G$, and  since $G$ is vertex-transitive the 
corresponding block system $\B_3$ is trivial with the whole set of vertices as the only block
(and $\Delta_3$ is universal). 
Consider $W_{-3}=\langle (1,-1,1,-1,1,-1,1,-1)\rangle$. It follows that
$N_{-3}$ is the largest subgroup in $G$  preserving the sets 
$\{b_1,b_2,b_3,b_4\}$ and $\{b_1^1,b_2^1,b_3^1,b_4^1\}$ 
(here $b_i^1$ denotes the image of the vertex $b_i$ under $1 \in H = \ZZ_{2}$);
these two sets coincide with the bipartition sets of $Q_3$. 
Note that the group $N_{-3}$ 
has to be of index $2$ in $G$, and therefore acts transitively
on each of the two sets of bipartition. Hence,
the corresponding block system $\B_{-3}$ coincides with the bipartition of $Q_3$
(and $\Delta_{-3}$ is the universal partition).
}
\end{example}

%%%%%%%%%%%%%%%%%%%%%%%%%%%
%%
%%      SECTION 5
%%
%%%%%%%%%%%%%%%%%%%%%%%%%%%

\section{Existence of imprimitivity block systems}
\label{sec:five}
\indent

Recall from Proposition \ref{cor:spec} that the spectrum of a $G$-invariant
digraph $\Ga$ is the union of spectra of all matrices $\chi(\S)$, $\chi \in H^*$, where
$\S$ is a symbol of $\Ga$. Recall that the matrix $\chi_0(\S)$ is nothing but the
adjacency matrix of $\Ga / H$. Also, the condition on the multiplicities of $\lm$ as in
Theorem~\ref{thm:Extreme} is equivalent to saying that $K_{\S, \lm} = \{ \chi_0\}$.
This condition is further equivalent to
saying that $H \leq N_\lm$, where  $N_\lm$ is the kernel of  the action of $G$ on the $\lm$-eigenspace
$W_{\lm}$ (see Lemma~\ref{lem:K}(i)) --
which is in turn equal to the kernel of the action of $G$ on the orbits of $N_\lm$
(see Lemma~\ref{lem:K}(ii)). Instead of Theorem~\ref{thm:Extreme}  we now prove
the following more general result.

\begin{theorem}
\label{thm:moregeneral}
Let $G$ be a transitive permutation group on a finite set $X$ containing an abelian
semi\-regular subgroup $H$. Further, let $\Ga$ be a $G$-invariant digraph with symbol $\S$,
and let $\lm \neq \val(\Ga)$ be an eigenvalue of $\Ga$. Let $N_\lm$ be the kernel of the action
of $G$ on the $\lm$-eigenspace $W_\lm$.
Then
\begin{enumerate}
\item[(i)] $N_\lm$ contains $\langle K_{\S,\lm} \rangle^\bot$.
\item[(ii)] If $\la K_{\S,\lm} \ra < H^*$ is a proper subgroup, then the orbits of $N_\lm$ form a
nontrivial $G$-invariant partition of $X$.
\end{enumerate}
\end{theorem}

\proof
Part (i) follows directly from Lemma \ref{lem:K}(i). For part (ii),
let $\B_\lm$ denote the $G$-invariant partition of $X$ consisting of
the orbits of $N_\lm$. By Proposition~\ref{lem:gtriple},
$$%
\B_\lm=\Pi(\un{x},\Delta_\lm,\phi_\lm)
$$%\
for an appropriate choice of $\un{x}$, $\Delta_\lm$, and $\phi_\lm$.
Since  $\bot \colon \Sub(H) \to \Sub(H^*)$  is an involutory antiisomorphism
between the subgroup lattices  of $H$ and $H^*$, we have that
$ \la K_{\S,\lm} \ra^\bot$ is nontrivial.
Hence $N_\lm$ is nontrivial by (i), implying that
$\B_\lm$ does not consist of singletons.

We  now show that $N_\lm$ is intransitive.
This is clear if $N_\lm \cap H < H$.
Suppose that  $N_\lm \cap H = H$ (which is equivalent to  $K_{\S,\lm} =\{ \chi_0\}$, by Lemma \ref{lem:K}(i)).
By Corollary~\ref{cor:gtriple} we have that  $\phi_{\lm}(i)=H$ for all $i$. 
Since $\lm \neq \val(\Ga)$  it follows that
$(1,1, \ldots, 1) \not \in V_{\lm, \chi_0}$. Thus, from the definition of
$\Delta_{\lm,\chi_0}$ it follows that
$\Delta_{\lm,\chi_0} \neq \Delta_U$. Hence $\Delta_{\lm} \neq \Delta_U$,  by Lemma~\ref{lem:Delta},
and  $N_\lm$ is intransitive.
\Qed

A transitive permutation group $G$ is called {\em quasiprimitive} if every nontrivial normal 
subgroup of $G$ is transitive.
An immediate consequence of Theorem~\ref{thm:moregeneral} is the
following result dealing with quasi\-pri\-mitive groups containing a semiregular abelian subgroup.
A special case of Corollary \ref{thm:primitive}, where  $G$ is primitive and $H$ is regular,
was proved by Knapp \cite{K}.

\begin{corollary}
\label{thm:primitive}
Let $G$ be a quasiprimitive permutation group containing an abelian semi\-regular
subgroup $H$, and let $\Ga$ be a $G$-invariant digraph with symbol $\S$.  If
$\lm \neq \val(\Ga)$ is  an eigenvalue of $\Ga$, then  $\la K_{\S,\lm} \ra = H^*$.
\end{corollary}

\bigskip
 Generalizing the
notion of a mixer of a bicirculant introduced in Section~\ref{sec:one},
we say that $g \in G$ is a
{\em mixer relative to} $H$ (in short, a {\em mixer} when the subgroup $H$ is
clear from the context), if the orbits of $H$ are not blocks of imprimitivity for
$\la g \ra$. The following corollary is just a rephrasing
of a special case covered in Theorem~\ref{thm:moregeneral}
into the mixers language.

\begin{corollary}\label{cor:mix}
Let $G$ be a transitive permutation group on a finite set $X$
containing an abelian  semiregular subgroup $H$ with a prime number of orbits, and let
$\Gamma$ be a $G$-invariant digraph with symbol $\S$.
If $G$ has a mixer relative to $H$, then for any eigenvalue $\lm \neq \val(\Ga)$
such that $\chi_0 \in K_{\S,\lm}$
we have $K_{\S,\lm} \neq \{ \chi_0\}$.
\end{corollary}

\proof
Suppose, on the contrary, that $K_{\S,\lm} = \{ \chi_0\}$. Then $N_\lm$ is intransitive,
by Theorem~\ref{thm:moregeneral}(ii), and $N_\lm \geq \la \chi_0 \ra^{\bot} = H$,  
by  Theorem~\ref{thm:moregeneral}(i). Therefore $\phi(i) = H$ for all $i$, and since $N_\lm$ is intransitive,
$\Delta_\lm \neq \Delta_U$. As $H$ has prime number of orbits, it follows that $\Delta_\lm = \Delta_I$,
and $G$ has no mixers relative to $H$, a contradiction. \Qed

\bigskip
When the number of orbits of $H$ is relatively small,
some additional information about mixers can be
obtained using a more direct combinatorial approach which takes into account  the intersections of orbits
of $H$ with their images under a mixer.
We here consider the case when $H$ has two orbits.
Let $\Ga$ be  a $G$-invariant digraph
(in this context also known as a {\em bi-Cayley digraph}) with symbol $\S$.
For convenience we  use abbreviations  $S = S_{11}$, $T = S_{12}$, $Q = S_{21}$ and $R = S_{22}$.
Note that $|S| = |R|$ and $|T| = |Q|$. Therefore, the eigenvalues of $\chi_0(\S)$ are
$\val(\Ga) = |S|+|T|$ and $d(\Ga)=|S| - |T|$.
For $g \in G$ we define the subsets $L_g$ and $M_g$ of $H$ by
\begin{eqnarray*}
L_g &=& \{h \in H \mid (b_1^h)^g \in X_1  \}, \\
M_g &=& \{h \in H \mid (b_2^h)^g \in X_1  \}.
\end{eqnarray*}
Since $g$ is an automorphism we have  $|L_g|+|M_g|=|H|$. Observe that  $g$ is a swap if and only if
$L_g=\emptyset$ and $M_g=H$, and that
$g$ is a mixer if and only if $L_g \neq \emptyset$ and $L_g \neq H$.

In what follows it will be convenient to view subsets of $H$ as elements of
the group algebra $\QQ H$. Following \cite{W}, a subset $A$ of $H$ is considered as
$\un{A} = \sum_{h\in A}h \in \QQ H$.
Characters of $H$ are naturally extended to an algebra-homomorphism $\QQ H \to \CC$ by letting
$\chi(\alpha)=\sum_{h\in H}a_h\chi(h)$ for $\alpha=\sum_{h\in H}a_h h$. 
Note  that, given a subset $S \subseteq H,$ we have $\chi(S) = \sum_{g\in S}\chi(g^{-1}) = \chi(\un{S^{-1}})$.
\medskip

\begin{theorem}\label{thm:MN}
Let $G$ be a transitive permutation group containing an abelian
semiregular subgroup $H$ having two orbits,
and let $\Ga$ be a $G$-invariant  digraph.
Then with notation and assumptions above, for all $g \in G$ and
for all characters $\chi \in H^*\setminus K_{\S,d(\Ga)}$ we have
$$
\chi(L_g)=\chi(M_g)=0.
$$
\end{theorem}

\proof
Let $g \in G$. For convenience we use a shorthand notation $L=L_g$ and $M=M_g$.
Choose  $h\in H$. The number of out-neighbours of $b_1^{h}$ that are contained in  $b_1^L \cup b_1^M$
is equal to $|hS \cap L| + | hT \cap M|$. This number is equal to
the coefficient of $h$ in $\un{S^{-1}} \;\un{L} + \un{T^{-1}} \;\un{M} \in \QQ H$.
Moreover, observe that it  is also equal to the number of out-neighbours of
$b_1^{hg}$ that are contained in $X_1$, which is either $|S|$ or $|Q|$ ($= |T|$), depending on
whether $h \in L$ or $h \in H\setminus L$. Consequently,
$$%
\un{S^{-1}} \; \un{L} + \un{T^{-1}} \; \un{M}=|S| \; \un{L} + |T| \;\un{H \setminus L}=
 d(\Ga) \; \un{L} + |T| \; \un{H}.
$$%
Similarly, by counting the number of out-neighbours of $b_2^{h}$ that are contained in $b_1^L \cup b_1^M$
(which equals the number of out-neighbours of $b_2^{hg}$ that are contained in $X_1$) we obtain
the equation
$$%
\un{Q^{-1}} \; \un{L} + \un{R^{-1}} \; \un{M}=|S| \; \un{M} + |T| \;\un{H \setminus M}=
 d(\Ga) \; \un{M} + |T| \; \un{H}.
$$%
Applying  a nonprincipal character $\chi \in H^*$ to both of the above
equalities we find  that $\chi(\un{L})$ and
$\chi(\un{M})$ are solutions  of the following   system of linear equations

$$
\begin{array}{rcrccc}
(\chi(\un{S^{-1}})- d(\Ga)) & x \, +& \chi(\un{T^{-1}}) & y  & = & 0\; \\
\chi(\un{Q^{-1}}) & x \, +& (\chi(\un{R^{-1}})- d(\Ga)) &y & = & 0.
\end{array}
$$
If $\chi \not\in K_{\S, d(\Ga)}$ we have
$\Det(\chi(\S) - d(\Ga) I) = (\chi(S)- d(\Ga)) (\chi(R)- d(\Ga)) - \chi(T)\chi(Q) \neq 0$.
As  the determinant of the above system 
is equal to %the conjugate of 
$\Det(\chi(\S) - d(\Ga) I)$, it follows that $\chi(L) = \chi(M) = 0$.
\Qed

\bigskip
We remark that in the particular case  when $H$ has  two orbits,
Corollary~\ref{cor:mix} follows directly from Theorem~\ref{thm:MN}. Namely,
let  $g \in G$ be a mixer relative to $H$. Then $L_g$ and $M_g$ are nonempty proper subsets of $H$.
It is well known that if all non-principal characters vanish on a subset $A \subseteq H$, then
either $A = \emptyset$ or $A = H$; see for instance \cite[Corollary~1.3.5]{BS}. Therefore,
there exists a non-principal character $\chi \in H^*$ such that either
$\chi(L_g) \neq 0$ or $\chi(M_g) \neq 0$. By Theorem~\ref{thm:MN}, we have
$\chi \in K_{\S, d(\Ga)}$.

\bigskip

%%%%%%%%%%%%%%%%%%%%%%%%%%%%%
%%
%%        SECTION 6
%%
%%%%%%%%%%%%%%%%%%%%%%%%%%%%%

\section{An application: generalized Petersen graphs}
\label{sec:seven}

\indent
To further illustrate the possible use of  techniques developed here  we
give an alternative short proof
of a well know classical result about edge-transitivity
(and thus arc-transitivity) of generalized Petersen graphs
due to Frucht, Graver and Watkins \cite[Theorem~2]{FGW}.
Recall that the {\em generalized Petersen graph} $\GP(n,s)$ is
the bicirculant relative to  the cyclic group $H = \ZZ_n$ and
admitting a symbol $\S$ with
$S_{11}=\{1,-1\}$, $S_{12}=S_{21}=\{0\}$ and $S_{22}=\{ s,-s\}$
where $s \in \ZZ_n\setminus \{0, n/2\}$.

\begin{theorem}
\label{t:gp}
A generalized Petersen graph $\Ga=\GP(n,s)$ is
edge-transitive (and thus arc-transitive)
if and only if $(n,s)$ is one of the following pairs:
$(4,1)$, $(5,2)$, $(8,3)$, $(10,2)$, $(10,3)$, $(12,5)$, or $(24,5)$.
\end{theorem}

\def\wt{\widetilde}

As remarked in Section~\ref{sec:one}, theorems of this  nature are  usually proved  either by
an elementary, although technically rather involved, combinatorial  approach 
(as, for instance, counting the number of $8$-cycles in the original proof of
the above theorem), or else by using  a two-step (normal subgroup) reduction  approach.
The first step identifies a  small restricted  family  from which all other graphs in the
class can be reconstructed. The actual constructions are then provided in the second step. 
In our proof we will make use of graph covers, and therefore we briefly review  
the necessary definitions and facts. \medskip

A {\em covering projection} of a graph $\wt{X}$ is a surjective mapping $p \colon \wt{X}  \to X$
such that for each $\tilde{u} \in V(\wt{X} )$ the set of arcs emanating from $\tilde{u}$ is
mapped bijectively onto the set of arcs emanating from $u = p(\tilde{u})$. The graph $\wt{X}$
is called a {\em covering graph} of the {\em base graph} $X$. The set 
$\textrm{fib}_u =p^{-1}(u)$ is a {\em fibre} of a vertex $u \in V(X)$. The subgroup $K$ of all those
automorphisms of $\tilde{X}$ which fix each of the fibres setwise is called the {\em  group
of covering transformations}. The graph $\tilde{X}$ is also called a {\em $K$-cover} of
$X$. It is a simple observation that the group of covering transformations of a connected
covering graph acts semiregularly on each of the fibres. In particular, if the group of
covering transformations is regular on the fibres of $\tilde{X}$, we say that $\tilde{X}$
is a {\em regular $K$-cover}. We say that $\alpha\in \Aut X$ {\em lifts} to an automorphism
of $\tilde{X}$ if there exists an automorphism $\tilde{\alpha}\in\Aut \tilde{X}$, called
a {\em lift} of $\alpha$, such that $\tilde{\alpha}p=p\alpha$. If the covering graph $\wt{X}$
is connected then  $K$ is the lift of the trivial subgroup of $\Aut X$. Note that a subgroup
$G \le  \Aut \wt{X}$ projects if and only if the partition of $V(\tilde{X})$ into the
orbits of $K$ is $G$-invariant.

A combinatorial description of a $K$-cover was introduced through a voltage graph by Gross and
Tucker \cite{GrT} as follows. Let $X$ be a graph and $K$ be a finite group. A {\em voltage
assignment} of $X$ is a mapping $\zeta \colon A(X) \to K$
with the property that $\zeta(u,v)=\zeta(v,u)^{-1}$ for any arc $(u,v) \in A(X)$   
(here, and in the rest of the paper,  $\zeta(u,v)$ is written instead of $\zeta((u,v))$ for the sake of brevity). 
The voltage assignment $\zeta$  extends to walks in $X$ in a natural way. 
In particular, for any walk $D=u_1u_2 \cdots u_t$ of $X$
we let $\zeta(D)$ denote the product voltage $\zeta({u_1,u_2})\zeta({u_2,u_3})\cdots \zeta({u_{t-1},u_{t}})$ of 
$D$, that is, the $\zeta$-voltage of $D$. 
The values of $\zeta$ are called {\em voltages}, and $K$ is the {\em voltage group}. The {\em voltage
graph} $X \times_{\zeta} K$ derived from a voltage assignment $\zeta \colon A(X) \to K$ has
vertex set $V(X) \times K$, and edges of the form $ \{(u,g),(v,\zeta(x)g)\}$, where
$x=(u,v) \in A(X)$. Clearly, $X \times_\zeta K$ is a covering of $X$ with the first
coordinate projection. By letting $K$ act on $V(X \times_\zeta K)$ as $(u,g)^{g'}=(u,gg')$,
$(u,g) \in V(X \times_\zeta K)$, $g' \in K$, one obtains a semiregular group of automorphisms
of $X \times_\zeta K$, showing that $X \times_\zeta K$ can in fact be viewed as a $K$-cover
of $X$.
 
Given a spanning tree $T$ of $X$, the voltage assignment $\zeta\colon A(X)\to K$ is said
to be {\em $T$-reduced} if the voltages on the tree arcs equal  the identity element in $K$.
In \cite{GrT2} it is shown that every regular covering graph $\tilde{X}$ of a graph $X$
can be derived from a $T$-reduced voltage assignment $\zeta$ with respect to an arbitrary
fixed spanning tree $T$ of $X$. 

The problem of whether an automorphism $\alpha$ of $X$ lifts
or not is expressed in terms of voltages as follows. 
Given $\alpha\in\Aut X$
and the set of fundamental closed  walks $\mathcal{C}$ based at a fixed
vertex $v\in V (X)$, we define  
$\bar{\alpha}=\{(\zeta(C),\zeta(C^\alpha))\mid C\in\mathcal{C}\}\subseteq K\times K$. 
Note that if $K$ is
abelian, $\bar{\alpha}$ does not depend on the choice of the base vertex, and the fundamental
closed walks at $v$ can be substituted by the fundamental cycles generated by the cotree arcs
of $X$. Also, from the definition, it is clear that for  a $T$-reduced voltage assignment
$\zeta$ the derived graph $X\times_\zeta K$ is connected if and only if the voltages of the cotree
arcs generate the voltage group $K$. It was proved in \cite[Theorem~4.2]{M} that, 
given  a connected regular cover $X \times_\zeta K$  of a graph $X$
derived from a voltage assignment $\zeta$  with the voltage group $K,$ then an automorphism
$\alpha$ of $X$ lifts if and only if $\bar{\alpha}$ is a function which extends to an automorphism of $K$.

\bigskip
\noindent
{\bf Proof of Theorem~6.1.} \quad
Our proof   uses the above two steps approach. In the first part we
show that every edge-transitive generalized Petersen graph $\Ga$
is a regular cyclic cover either of the cube or of the
Petersen graph (such that $\Aut(\Ga)$ projects). In the second part we
determine all such covers.

\medskip
{\sc Part 1.}
First  note that the eigenvalues of $\chi_0(\S)$ are
$\val(\Ga) = 3$ and $d(\Ga) = 1$.
The respective eigenspace $V_{1,\chi_0}$ is  spanned by $(1,-1)$.
By the definition of the relation $\sim_{1,\chi_0}$ the partition
$\Delta_{1,\chi_0}$  is equal to $\Delta_I$,
and by Lemma~\ref{lem:Delta}, the partition
$\Delta_1$ is also equal to $\Delta_I$.
Let $K=N_1 \cap H$. By Lemma~\ref{lem:K}(iii), $\phi(1)=\phi(2)=K$.
Thus, the $\Aut(\Ga)$-invariant partition $\B_1$ of the
vertex set of $\Ga$, induced by the eigenvalue $1$, is just the partition formed by
the orbits of $K$. Since the two orbits of $H$ are joined by one perfect
matching, we get that any two orbits of $K$ are joined by at most one perfect
matching. This implies that $N_1 = K$ and  $\Gamma$ is a normal cover of
$\Gamma/K$. 

\medskip
In order to determine $K$, and consequently $\Ga/K$,
recall that
$K = \la K_{\S, 1}\ra^\bot = \cap \Ker \chi$,
where $\chi$  ranges over all characters in $K_{\S, 1}$.
Therefore, we first find $K_{\S, 1}$.
The set $K_{\S, 1}$ is determined using the assumption that $\Ga$ is
edge-transitive, and hence that it has a mixer relative
to $H$. By  Corollary~\ref{cor:mix}, the set $K_{\S, 1}$ contains
a non-principal character, say $\chi $, implying that
\begin{equation}
\label{eq:char0}
(\chi(1)+\chi(-1)-1)(\chi(s)+\chi(-s)-1)=1.
\end{equation}
All solutions of (\ref{eq:char0})  give all possible characters in $K_{\S, 1}$. (Observe that
$\chi \in K_{\S, 1}$ if and only if $\chi^{-1} \in K_{\S, 1}$.)
For a natural number $m$ let  $\xi_m = e^{2\pi \im \over m}$.
As $\chi(1) = \xi_n^j$ for some
$j \in \ZZ_n\setminus \{0\}$,   equation  (\ref{eq:char0})   rewrites as
$(\xi_n^j+ \xi_n^{-j} - 1)(\xi_n^{js}+ \xi_n^{-js}-1)=1$.
By setting  $j' = j/\gcd(j,n)$, $n' = n/\gcd(j,n)$, and $s' = s\, \mod\,n'$,  this last  equation becomes
$$
(\xi_{n'}^{j'}+ \xi_{n'}^{-j'} - 1)(\xi_{n'}^{j's'}+ \xi_{n'}^{-j's'}-1)=1.
$$
Since $j'$ and $n'$ are coprime,
$\xi_{n'}^{j'}$ and $\xi_{n'}$ are both primitive $n'$-th roots of unity.
Hence their minimal polynomial is the same,
implying that
\begin{equation}
\label{eq:basiceq}
(\xi_{n'}+ \xi_{n'}^{-1} - 1)(\xi_{n'}^{s'}+ \xi_{n'}^{-s'}-1)=1.
\end{equation}
If both factors in  (\ref{eq:basiceq})  are positive then $n'=1$,  forcing $j = n$.
However, this is not possible as
$\chi$ is non-principal. Therefore  $\xi_{n'}+ \xi_{n'}^{-1} - 1 = 2\cos(2\pi/n')  - 1 <0$,  and so
$n' \leq 5$. A direct analysis yields that either $n'=4$  and $s' = \pm 1$, or else  $n'=5$
and $s' = \pm 2$. It follows that $j'= \pm 1$ in the first case
while $j' \in \{\pm 1, \pm 2\}$ in the second case.
Hence for some natural number $m$ either  $n = 4m$, $j = \pm m$,  and $s = \pm 1 (\mod\,4)$,  or else
$n = 5m$, $j \in \{\pm m, \pm 2m\}$, and $s = \pm 2 (\mod\,5)$.
 We now analyze each of these cases separately.

\medskip
Let $n = 4m$, $j = \pm m$, and $s = \pm 1 (\mod\,4)$, where $m$ is not divisible by $5$. Since
$\chi(1) = \xi_n^j \in \{\xi_4, \xi_4^{-1}\}$ we have that
$K_{\S, 1} = \{\chi_0, \chi, \chi^{-1}\}$, and so $K = \la K_{\S, 1} \ra^{\bot}$ is the
unique  index $4$ subgroup  in  $\ZZ_n$.   Since  $s = \pm 1 (\mod\,4)$,  the graph
$\Ga/K$ is the cube.

Let $n = 5m$, $j \in \{\pm m, \pm 2m\}$,  and $s = \pm 2 (\mod\,5)$ where $m$ is not divisible by $4$.
Since $\chi(1) = \xi_n^j \in \{\xi_5, \xi_5^{-1},  \xi_5^2, \xi_5^{-2}\}$ we have that
$K_{\S, 1} = \{\chi_0, \chi, \chi^{-1}, \chi^2, \chi^{-2}\}$,
and so $K = \la K_{\S, 1} \ra^{\bot}$ is the
unique  index $5$ subgroup  in  $\ZZ_n$.  Since $s = \pm 2 (\mod\,5)$,  the graph
$\Ga/K$ is the Petersen graph.

Let $n = 20 m$  where $s = \pm 1 (\mod\,4)$ and  $s = \pm 2 (\mod\,5)$.   By the above analysis,
the group $\la K_{\S, 1} \ra = \la \chi, \chi' \ra$ where $\chi(1) = \xi_4$ and
$\chi'(1) = \xi_5$. In this case $K = \la K_{\S, 1} \ra^{\bot}$ is the unique index $20$
subgroup in $\ZZ_n$. Since $s = \pm 1 (\mod\,4)$
and  $s = \pm 2 (\mod\,5)$ we have, by the Chinese
Remainder Theorem, that  $s = \pm 3\, (\mod\, 20)$ or $s = \pm 7\, (\mod\, 20)$.
Hence $\Ga/K$ is either $\GP(20, 3)$ or $\GP(20,7)$. In particular, $\Ga/K$ and hence also
$\Ga$ is a regular cyclic cover of both the cube and the Petersen graph.
(Actually, as it will become clear in Part 2, this case
cannot occur.)

Note that in all of the above cases $K$ is isomorphic to  $\ZZ_m$.

\medskip
{\sc Part 2.}
Recall from Part 1 that $\Ga \to \Ga/K$ is a regular cyclic
covering projection such that $\Aut(\Ga)$ projects,
and hence some  minimal arc transitive subgroup of $\Ga/K$ lifts. We determine all such covers
where the  base graph is either the cube or the Petersen graph. Since $\GP(20, 3)$ or
$\GP(20,7)$ are not among them, the case where the base graph is one of these two need
not be considered.

Let $Y$ be the cube. Denote its  vertices by
$\{0,1,2,3,4,5,6,7\}$ in such a way that $01230$   and  $45674$ are
the outer and the inner cycles, respectively, and $04$, $15$, $26$ and $37$ are the spokes. Choose
a  spanning tree consisting of edges  $01$, $12$, $23$, $04$, $15$, $26$,  and $37$.
The base homology cycles are
$C_1 = 30123$,
$C_2 = 45104$,
$C_3 = 56215$,
$C_4 = 67326$, and
$C_5 = 7401237$.
Voltages on the spanning tree are $0$.
Without  loss of generality we can assume that
the voltages on the remaining arcs are $\zeta(3,0) = 1$, and
$\zeta(4,5) = \zeta(5,6) = \zeta(6,7) = a$, and $\zeta(7,4) = a+1$.
(This is obtained using the fact that
 the derived covering graph is the generalized Petersen graph.)
Consider now  the automorphism $\alpha = (134)(527)$
which must lift since an edge-transitive group of $\Ga$ projects.
Thus by \cite[Theorem~4.2]{M}, the   
induced mapping $\bar{\alpha}$ must extend to an automorphism $x \mapsto \lm x$ of $K$ (see the above discussion). 
As $\bar{\alpha}$ maps $\zeta(C_j)$ to $\zeta(C_j^\alpha),$ where $j \in \{1,\ldots,5\},$ 
we get the following system of equations in $K$:
$$%
\lambda  =  a,  \quad
\lambda a  = 1, \quad
\lambda a  = a, \quad
\lambda a   = -4a -1, \quad
\lambda (a+1) = 2a.
$$%
This gives $\lambda = a = 1$ and $6 = 0$ in $K$,
implying that $K$ is isomorphic to $\ZZ_2$, $\ZZ_3$, or $\ZZ_6$.
If $K = \ZZ_2$ then $\Ga \cong \GP(8,5) \cong \GP(8,3)$. If $K = \ZZ_3$ then
$\Ga \cong \GP(12,5)$.
If $K = \ZZ_6$ then $\Ga \cong \GP(24,5)$.

Let $Y$ be the Petersen graph. Denote its  vertices by
$\{0,1,2,3,4,5,6,7,8,9\}$ in such a way that $012340$   and  $567895$ are
the outer and the inner cycles, respectively, and $05$, $16$, $27$, $38$,
and $49$ are the spokes. Choose
a  spanning tree consisting of edges
$01$, $12$, $23$, $34$, $05$, $16$, $27$, $38$, and $49$. The base homology cycles are
$C_1 = 401234$,
$C_2 = 572105$,
$C_3 = 683216$,
$C_4 = 794327$,
$C_5  =8501238$, and
$C_6 = 9612349$.
Voltages on the spanning tree are $0$.
Without  loss of generality we can assume that the voltages on the remaining arcs are $\zeta(4,0) = 1$, and
$\zeta(5,7) = \zeta(6,8) = \zeta(7,9) = a$, and $\zeta(8,5) = \zeta(9,6) = a+1$.
(This is obtained using the fact that
 the derived covering graph is the generalized Petersen graph.)
Consider now  the automorphism $\alpha = (154)(289)(367)$
which must lift since an edge-transitive group of $\Ga$ projects.
Since the induced mapping $\bar{\alpha}$ must extend to an automorphism $x \mapsto \lm x$ of $L,$ 
and also  $\bar{\alpha}$ maps $\zeta(C_j)$ to $\zeta(C_j^\alpha),$ where $j \in \{1,\ldots,6\},$ 
we get the following system of equations in $L$:
$$%
\lambda  =  -2a - 1,  \,
\lambda a  = a, \,
\lambda a  = 5a + 2, \,
\lambda a   = a, \,
\lambda (a+1) = -3a-1, \,
\lambda (a+1) = -3a-1.
$$%
This gives $2 = 0$, $\lambda = -1$ and either $a = 0$ or $a = 1$
in $K$, implying that $K$ is isomorphic to $\ZZ_2$.
If $a = 0$ then $\Ga \cong \GP(10,2)$. If $a = 1$ then $\Ga \cong \GP(10,7) \cong \GP(10,3)$,
completing the proof of  Theorem~6.1.
\Qed

\begin{footnotesize}
%\begin{small}

\end{footnotesize}
\end{document}